\documentclass[12pt]{amsart}
\usepackage{enumerate}
\theoremstyle{plain}

\newtheorem{thm}{Theorem}[section]
\newtheorem{cor}[thm]{Corollary}
\newtheorem{lem}[thm]{Lemma}
\newtheorem{prop}[thm]{Proposition}

\theoremstyle{definition}
\newtheorem{defi}[thm]{Definition}
\newtheorem{defis}[thm]{Definitions}
\newtheorem{conj}[thm]{Problem}
\newtheorem{conv}[thm]{Convention}
\newtheorem{nota}[thm]{Notation}
\newtheorem{rem}[thm]{Remark}
\newtheorem{rems}[thm]{Remarks}
\newtheorem{exa}[thm]{Example}
\newtheorem{exas}[thm]{Examples}
\newtheorem{sit}[thm]{}

\newcommand{\brem}{\begin{rem}}
\newcommand{\brems}{\begin{rems}}
\newcommand{\erem}{\end{rem}}
\newcommand{\erems}{\end{rems}}
\newcommand{\bexa}{\begin{exa}}
\newcommand{\bexas}{\begin{exas}}
\newcommand{\eexa}{\end{exa}}
\newcommand{\eexas}{\end{exas}}
\newcommand{\bdefi}{\begin{defi}}
\newcommand{\edefi}{\end{defi}}
\newcommand{\bdefis}{\begin{defis}}
\newcommand{\edefis}{\end{defis}}
\newcommand{\bcor}{\begin{cor}}
\newcommand{\ecor}{\end{cor}}
\newcommand{\blem}{\begin{lem}}
\newcommand{\elem}{\end{lem}}
\newcommand{\bconv}{\begin{conv}}
\newcommand{\econv}{\end{conv}}
\newcommand{\bconj}{\begin{conj}}
\newcommand{\econj}{\end{conj}}
\newcommand{\bprop}{\begin{prop}}
\newcommand{\eprop}{\end{prop}}
\newcommand{\bthm}{\begin{thm}}
\newcommand{\ethm}{\end{thm}}
\newcommand{\bnota}{\begin{nota}}
\newcommand{\enota}{\end{nota}}
\newcommand{\bsit}{\begin{sit}}
\newcommand{\esit}{\end{sit}}
\newcommand{\be}{\begin{equation}}
\newcommand{\ee}{\end{equation}}
\newcommand{\bproof}{\begin{proof}}
\newcommand{\eproof}{\end{proof}}
\def\ba{\begin{array}}
\def\ea{\end{array}}
\def\bea{\begin{eqnarray}}
\def\eea{\end{eqnarray}}

\def\bnum{\begin{enumerate}}
\def\enum{\end{enumerate}}

\newcommand{\HOT}{\operatorname{HOT}}

\newcommand{\Sy}{{\operatorname{Symm}}}

\newcommand{\Poly}{{\operatorname{Poly}}}

\newcommand{\A}{{\mathbb A}}

\newcommand{\PP}{{\mathbb P}}

\newcommand{\C}{{\mathbb C}}

\newcommand{\N}{{\mathbb N}}

\newcommand{\F}{{\mathbb F}}

\newcommand{\Mor}{{\operatorname{Mor}}}

\newcommand{\Hol}{{\operatorname{Hol}}}
\newcommand{\codim}{{\operatorname{codim}}}

\renewcommand{\phi}{\varphi}

\title{Hyperbolicity of general deformations: proofs}

\author{Mikhail Zaidenberg}
\address{Universit{\'e}
Grenoble I, Institut Fourier, UMR 5582 CNRS-UJF, BP 74, 38402 St.\
Martin d'H{\`e}res c{\'e}dex, France}
\email{zaidenbe@ujf-grenoble.fr}

\thanks{
\mbox{\hspace{11pt}}{\it 2000 Mathematics Subject Classification}:
14J70,  32J25.\\
\mbox{\hspace{11pt}}{\it Key words}: Kobayashi hyperbolicity,
projective hypersurface, deformation}
\thanks{
{\bf Acknowledgement:} This paper was written during a visit of
the author the Max-Planck-Institute of Mathematics, Bonn. He thank
this institution for a generous support and excellent working
conditions.}
\date{}

\begin{document}
\begin{abstract}
We  modify the deformation method from \cite{SZ} in order to
construct further examples of Kobayashi hyperbolic surfaces in
$\PP^3$ of any even degree $d\ge 8$.
\end{abstract}

\maketitle

Given a  hypersurface $X_d=f_d^*(0)$ in $\PP^n$ of degree $d$, we
say that a (very) general small deformation of $X_d$ is hyperbolic
if  for any (very) general degree $d$ hypersurface
$X_\infty=g_d^*(0)$ and for all sufficiently small
$\varepsilon\in\C\setminus \{0\}$ (depending on $X_\infty$) the
hypersurface $X_{d,\varepsilon}=(f_d+\varepsilon g_d)^*(0)$ is
Kobayashi hyperbolic. With this definition let us formulate the
following version of the Kobayashi Conjecture.

\smallskip

{\bf Weak Kobayashi Conjecture.} {\it For every hypersurface $X_d$
in $\PP^n$ of degree $d\ge 2n-1$, a (very) general small
deformation of $X_d$ is Kobayashi hyperbolic.}

\smallskip

The original Kobayashi Conjecture claims, in particular, that a
(very) general surface $X_d$ of degree $d\ge 5$ in $\PP^3$ is
Kobayashi hyperbolic. This is known to hold indeed for a very
general surface of degree at least $21$ (see McQuillan \cite{MQ}
and Demailly-El Goul \cite{DEG}).

By Brody's Theorem, a compact complex space $X$ is hyperbolic if
and only if any holomorphic map $\C\to X$ is constant. Hence the
proof of hyperbolicity reduces to a certain degeneration principle
for entire curves in $X$. The Green-Griffiths' proof of Bloch's
Conjecture \cite{GG} provides a kind of such degeneration
principle. According to this principle, every entire curve
$\varphi:\C\to X$ in a very general surface $X\subseteq\PP^3$ of
degree $d\ge 21$ satisfies an algebraic differential equation
\cite{DEG,MQ}. See also \cite{Ro, Si} for recent advances in
higher dimensions.

The deformation method showed to be quite effective to construct
examples of low degree hyperbolic surfaces in $\PP^3$. A nice
construction due to J.\ Duval \cite{Du} of a hyperbolic sextic
$X_{\varepsilon}\subseteq \PP^3$ uses this method iteratively in 5
steps, so that $\varepsilon=(\varepsilon_1,\ldots,\varepsilon_5)$
has 5 subsequently small enough components. Hence
$X_{\varepsilon}$ belongs to a 5-dimensional linear system;
however the deformation of $X_{0}$ to $X_{\varepsilon}$ neither is
linear nor very generic.

In \cite{SZ} we exhibited examples of some special surfaces $X_d$
in $\PP^3$ of any given degree $d\ge 8$ such that a {\it general}
small deformation of $X_d$ is Kobayashi hyperbolic. In these
examples $X_d=X_{d'}'\cup X_{d''}''$, where $d=d'+d''$, is a union
of two cones in $\PP^3$ with distinct vertices over plane
hyperbolic curves in general position.

Let us indicate briefly the deformation method used in \cite{SZ}
(see also the references in \cite{SZ,SZ1}). Given two
hypersurfaces $X_{d,0}$ and $X_{d,\infty}$  in $\PP^n$ of the same
degree $d$, we consider the pencil of hypersurfaces
$\{X_{d,\varepsilon}\}_{\varepsilon\in\C}$ generated by $X_{d,0}$
and $X_{d,\infty}$. Assuming that for a sequence $\varepsilon_n\to
0$, the hypersurfaces $X_{d,\varepsilon_n}$ are not hyperbolic,
there exists a sequence of Brody entire curves $\varphi_n:\C\to
X_{d,\varepsilon_n}$ which converges to a (non-constant) Brody
curve $\varphi:\C\to X_{d,0}$. Suppose in addition that the
hypersurface $X_{d,0}$ admits a rational map to a hyperbolic
variety $\pi:X_{d,0}\dashrightarrow Y_0$ (to a curve $Y_0$ of
genus $\ge 2$ in case where $\dim X_{d,0}=2$). Then necessarily
$\pi\circ \varphi=\,$cst, provided that the composition $\pi\circ
\varphi$ is well defined. Anyhow, the limiting Brody curve
$\varphi:\C\to X_{d,0}$ degenerates.

For a union $X_{d,0}=X_{d',0}'\cup X_{d'',0}''$ of two cones in
general position in $\PP^3$ as in \cite{SZ}, there is a further
degeneration principle. It prohibits to the image $\varphi(\C)$ to
meet the double curve $D=X_{d',0}'\cap X_{d'',0}''$ outside the
points of $D\cap X_{d,\infty}$. Using the assumptions that
$d',\,d''\ge 4$  and $X_{d,\infty}$ is general this forces
$\varphi$ to be constant, contrary to our construction.

This applies in particular to the union of two quartic cones
$X_{4,0}'\cup X_{4,0}''$ in $\PP^3$ in general position. Modifying
the construction in \cite{SZ}, in the present note we establish,
in particular, hyperbolicity of a general deformation of a double
quartic cone in $\PP^3$, see Example \ref{exa1} below.

 The author is grateful to the referee for indicating a flow in
the first draft of the paper.

 \section{Some technical lemmas}
 Here we expose some preliminary facts
 that will be used in the next section.
 We let $\Delta$ denote the unit disc in $\C$, $B^n$ the unit
 ball in $\C^n$ and $\Hol (B^n)$
 the space of all holomorphic functions on $B^n$. For
 two complex spaces $X$ and $Y$, $\Hol (X,Y)$
 stands for the space of
 all holomorphic maps $X\to Y$ with the usual topology.

 \blem\label{prop} Let $f_0,f_\infty\in \Hol (B^n)$
 be such that
  $f_0(0)=f_\infty(0)=0$ and the
 divisors $X_0=f_0^*(0)$ and $X_\infty=f_\infty^*(0)$
 have no common
 component passing through $0$.
 Let $\Gamma=X_0\cap X_\infty$ and
 $X_\varepsilon=f_\varepsilon^{-1}(0)$, where
 $f_\varepsilon=f_0+\varepsilon f_\infty$.
 We assume that $\bigtriangledown f_0 |_\Gamma=0$.
 Let further
 $\varphi_n\in\Hol (\Delta, X_{\varepsilon_n})$,
 where $\varepsilon_n
 \longrightarrow 0$, be a sequence
 of holomorphic discs which converges
 to $\varphi\in\Hol (\Delta, X_0)$
 with $\varphi(0)=0$. Then
 necessarily
 $d\varphi (0)\in T_0 X_\infty$.
 \elem

\bproof The assertion is clearly true in the case where
$\varphi(\Delta)\subseteq \Gamma$. So we will assume further that
$\varphi(\Delta)\not\subseteq \Gamma$.

\smallskip

\noindent {\it Claim 1. Under the assumptions as above
$\varphi_n(t_n)\in\Gamma$ for some sequence $t_n\longrightarrow
0$. }

\smallskip

\noindent {\it Proof of Claim 1.} Let us consider the holomorphic
map
$$F: B^n\to \C^2,\qquad z\longmapsto (f_0(z),f_\infty (z))\,.$$
It is easily seen that $F$ possesses the following properties:
\begin{enumerate}\item[$\bullet$]
$F(0)=0$; \item[$\bullet$] $F^{-1}(0)=\Gamma$; \item[$\bullet$]
$F(X_{\varepsilon_n})\subseteq l_n$, where $l_n:=\{x+\varepsilon_n
y=0\}\subseteq  \C^2$; \item[$\bullet$] $F(X_0)\subseteq
l_0:=\{x=0\}$; \item[$\bullet$] $F\circ \varphi_n(\Delta)\subseteq
l_n$; \item[$\bullet$] $F\circ \varphi(\Delta) \subseteq l_0$,
$F\circ \varphi (0)=0$, $F\circ \varphi\not\equiv 0$.
\end{enumerate}

\noindent We let $F\circ \varphi_n=(x_n(t), y_n(t))$ and $F\circ
\varphi=(0,y(t))$. Thus $x_n\longrightarrow 0$ and $
y_n\longrightarrow y$ as $n\longrightarrow \infty$. Since $y(0)=0$
and $y\not\equiv 0$, we have $y_n\not\equiv 0$. By Rouch\'e's
Theorem there exists a sequence $t_n\longrightarrow 0$ such that
$y_n(t_n)=0$, so also $x_n(t_n)=-\varepsilon_n y_n(t_n)=0$. Hence
$\varphi_n (t_n)\in\Gamma=X_0\cap X_\infty$, as claimed. \eproof

It will be convenient for the rest of the proof to replace the
given sequence $(\varphi_n)$ by a new one $(\psi_n)$. We let
$\psi_n(t)=\varphi_n(a_nt+t_n)$ with $(t_n)$ as in Claim 1 and
$a_n:=1-|t_n|\longrightarrow 1$. Then $\psi_n\in\Hol (\Delta,
X_{\varepsilon_n})$ and $\psi_n\longrightarrow \varphi$ as
$n\longrightarrow \infty$. Moreover
$p_n:=\psi_n(0)=\varphi_n(t_n)\in\Gamma$ $\forall n\ge 1$ and
$v_n:=d\psi_n(0)\longrightarrow v:=d\varphi(0)$ when
$n\longrightarrow \infty$. Now the assertion follows immediately
from the next claim.

\smallskip

\noindent {\it Claim 2. $v_n\in T_{p_n}X_\infty$ $\forall n\ge
1$.}

\smallskip

\noindent {\it Proof of Claim 2.} We have:
$$\psi_n(t)=p_n+tv_n+\HOT (t)\qquad\mbox{and}
\qquad f_{\varepsilon_n}(x) =\langle \bigtriangledown
f_{\varepsilon_n}(p_n), x-p_n\rangle + \HOT(x-p_n)\,,$$ where
$\HOT$ means ``the higher order terms". Hence \be\label{eq}
f_{\varepsilon_n}\circ\psi_n (t)=\langle \bigtriangledown
f_{\varepsilon_n}(p_n), v_n\rangle \cdot t + \HOT(t)\,.\ee Using
(\ref{eq}) and the identity $f_{\varepsilon_n}\circ\psi_n \equiv
0$ we obtain $$0=\langle \bigtriangledown f_{\varepsilon_n}(p_n),
v_n\rangle=\langle \bigtriangledown f_0(p_n), v_n\rangle +
\varepsilon_n\langle \bigtriangledown f_\infty(p_n),
v_n\rangle=\varepsilon_n\langle \bigtriangledown f_\infty(p_n),
v_n\rangle\,.$$ Indeed, by our assumption $\bigtriangledown
f_0|_\Gamma=0$, in particular $\bigtriangledown f_0(p_n)=0$
$\forall n\ge 1$. This proves the claim. \qed

\medskip
 Consider, for instance,  a pencil  of degree $d$
hypersurfaces
$$X_\varepsilon=(f_0+\varepsilon f_\infty)^*(0)
\qquad\mbox {in} \quad\PP^{n+1}$$ generated by
$$X_0=X_0'\cup X_0''=f_0^*(0)\quad\mbox{ and}
\quad X_\infty=f_\infty^*(0)\,.$$ Assume that $D:=X_0'\cap
X_0''\subseteq X_\infty$. Then for any sequence of entire curves
$\varphi_{n}:\C\to X_{\varepsilon_n}$ which converges to
$\varphi:\C\to X_0'$   we have by Lemma \ref{prop}:
$$d\varphi (t)\in T_P X_0'\cap T_P X_\infty\quad\forall P=\varphi
(t)\in D\,.$$

\smallskip

Next we study an enumeration problem, which deals with the
intersection of a general hypersurface and generators of a given
cone in $\PP^{n+1}$.

\bprop\label{trans} We let $\widehat Y\subseteq\PP^{n+1}$ be a
cone over a variety $Y\subseteq\PP^n$. We consider also a general
hypersurface $X\subseteq \PP^{n+1}$ of degree $e\ge 2\dim Y$. Then
$X$ meets every generator $l=(PQ)$ of $\widehat Y$, where $P$ is
the vertex of the cone and $Q$ runs over $Y$, in at least
$k=e-2\dim Y$ points transversally. \eprop

\bproof  We use below the following notation. For a pair
$(n,e)\in\N^2$ we let $\F(n+1,e)$ denote the vector space of all
homogeneous forms in $n+2$ variables of degree $e$ and
$\PP(n+1,e)$ its projectivization.  We let $CY$ denote the affine
cone over $Y$ and $CY^*=CY\setminus\{0\}$ the same cone with the
vertex deleted. Let us fix coordinates in $\PP^{n+1}$ in such a
way that $P=(0:\ldots:0:1)$ and $Y\subseteq \{z_{n+1}=0\}$. If
$Q=(z_0:\ldots :z_n:0)=(z':0)\in Y$ then
$$(PQ) = \{(z':z_{n+1})\,|\,z_{n+1}\in\C\}\cup
\{P\}\,.$$ For a hypersurface  $X$  in $\PP^{n+1}$ of degree $e$
its defining equation $f=0$ can be written in the form
\be\label{eq1} f(z',z_{n+1})=\sum_{i=0}^{e}
a_i(z')z_{n+1}^{e-i}=0\,, \ee where $a_i$ is a homogeneous form in
$z'$ of degree $i$. Assuming that $P\notin X$ i.e., $a_0\neq 0$,
we can normalize the equation so that $a_0=1$. Fixing
$z'\in\A^{n+1}$ we specialize $f$ to a monic polynomial
$f_{z'}\in\C[z_{n+1}]$ of degree $e$. In these terms the
proposition asserts that for $k=e-2\dim Y$ and for a general $f\in
\F(n+1,e)$, the specialization $f_{z'}$ has at least k simple
roots whatever is the choice of $z'\in CY^*\subseteq\A^{n+1}$.

The affine chart
$$U=\PP(n+1,e)\setminus\{a_0=0\}$$ can be identified
with the affine space of all sequences of homogeneous forms
$a=(a_1,\ldots,a_{e})$ with $\deg a_i=i$. The specialization
$(f,z')\longmapsto f_{z'}$ defines a morphism
$$\tilde{\rho}:U\times CY\to \Poly_{e}\,,$$
where $\Poly_{e}$ stands for the affine variety of all monic
polynomials of degree $e$. In turn $\Poly_{e}$ can be identified
with $\Sy_{e}(\A^1)\cong \A^{e}$.

Let us consider further the Vieta map $$\nu: \A^{e}\to
\Poly_{e},\qquad (\lambda_1,\ldots,\lambda_{e})\longmapsto p(z)
=\prod_{i=1}^{e} (z-\lambda_i)\,.$$ This is a ramified covering of
degree $e!$. For a multi-index $\bar n=(n_1,\ldots,n_s)$ with
$\sum_{i=1}^s n_i=e$ we let $$\Sigma'_{\bar n}=\nu (D_{\bar n})
\subseteq \Poly_{e}\,,$$ where $D_{\bar n}$ is the linear subspace
of $\A^{e}$ given by equations
$$\lambda_1=\ldots=\lambda_{n_1},\quad
\lambda_{n_1+1}=\ldots=\lambda_{n_1+n_2},
\quad\ldots,\quad\lambda_{n_1+\ldots+n_{s-1}+1}
=\ldots=\lambda_{e}\,.$$ Clearly both $D_{\bar n}$ and
$\Sigma'_{\bar n}$ have pure dimension $s$. Letting
$$\Sigma'_{k}=\bigcup_{n_k\ge 2}
\Sigma'_{\bar n}\subseteq \Poly_{e}\,$$ denote the variety of all
monic polynomials of degree $e$ with at most $k-1$ simple roots,
we have
$$\dim \Sigma'_k= \max_{n_k\ge 2}\,
\{\dim \Sigma'_{\bar n}\}= k-1+ \left[\frac{e-k+1}{2}\right]\,.$$
If $e-k+1$ is even then the latter maximum is achieved for
$$n_1=\ldots=n_{k-1}=1,\,\,\, n_k=\ldots =n_s=2\,,$$
and otherwise for
$$n_1=\ldots=n_{k-2}=1,\,\, \,n_{k-1}=\ldots
=n_s=2\,.$$ Anyhow
$$\codim\, (\Sigma'_k,\Poly_{e})
=1+\left[\frac{e-k}{2}\right]\,.$$

\smallskip

\noindent {\it Claim 1. The restriction $d\tilde\rho |_{TU}$ is
surjective at every point $(a,z')\in U\times CY^*$. In particular
$d\tilde\rho$ has maximal rank $e$ at every such point. }

\smallskip

\noindent {\it Proof of Claim 1.} For a point
$(a,z')=(a_1,\ldots,a_{e},z_0,\ldots,z_n) \in U\times CY^*$ we let
$$a^0=(a^0_1,\ldots,a^0_{e})\in\A^{e},\quad\mbox{ where}\quad
a_i^0=a_i(z'),\,\,i=1,\ldots,e\,.$$ Since $z'\neq 0$, for an
arbitrary tangent vector $b^0=(b^0_1,\ldots,b^0_{e})\in\A^{e}$
there exists a $e$-tuple of homogeneous forms
$b=(b_1,\ldots,b_{e})$ with $\deg b_i=i$ such that $b(z')=b^0$.
Therefore $$(a+tb)(z')=a^0+tb^0\quad\mbox{and so} \quad
d\tilde\rho (a^0,z')(b,0)=b^0\,.$$ This proves Claim 1. \qed

\smallskip

By virtue of Claim 1,
$$\codim\, (\tilde\rho^{-1}(\Sigma'_k),\,U\times CY^*)=
\codim\, (\Sigma'_k,\,\Poly_{e})
=1+\left[\frac{e-k}{2}\right]\,.$$ Since $$f_{\lambda
z'}(z_{n+1})=\lambda . f_{z'}(z_{n+1}) =\lambda^{-e}f_{z'}
(\lambda z_{n+1})\qquad\forall \lambda\in\C^*\,,$$ the subvariety
$\tilde\rho^{-1}(\Sigma'_k)$ of $U\times CY^*$ is stable under the
natural $\C^*$-action on the second factor. Hence
$$\codim\, \left(\tilde\rho^{-1}(\Sigma'_k)/\C^*,\,U\times Y\right)
=\codim\, (\tilde\rho^{-1}(\Sigma'_k),\,U\times CY^*)=
1+\left[\frac{e-k}{2}\right]\,.$$ Thus the general fibers of the
projection $${\rm pr}_2: U\times Y\to U\,$$ do not meet
$\tilde\rho^{-1} (\Sigma'_k)/\C^*\subseteq U\times Y$ provided
that
$$\dim Y\le\left[\frac{e-k}{2}\right]\,.$$
The latter inequality is equivalent to $k\le e-2\dim Y$, which
fits our assumption. Now the proposition follows. \eproof

\brem\label{gp} Let us indicate an alternative approach. Given a
projective variety $Y\subseteq\PP^n$ and a cone
$X\subseteq\PP^{n+1}$ over $Y$ with vertex $P$, for every $k\ge 1$
we consider the subset $\F(Y,e,k)\subseteq \F(n+1,e)$ of all forms
$f\in\F(n+1,e)$ such that the intersection divisor $f^*(0)\cdot
(PQ)$ has at most $k-1$ reduced points on at least one generator
$l=(PQ)$ ($Q\in Y$) of $X$. We let $\PP(Y,e,k)$ denote the
projectivization of $\F(Y,e,k)$. Proposition \ref{trans} asserts
that the complement $\PP(n+1,e)\setminus \PP(Y,e,k)$ is a nonempty
Zariski open subset of $\PP(n+1,e)$ provided that $e\ge 2\dim Y
+k$. We divide this into two claims; the first one is proved in a
general setting, while for the second one we provide a simple
argument in dimension 3 only.

\smallskip

{\it Claim 1. $\PP(Y,e,k)$ is a Zariski closed subset of
$\PP(n+1,e)$.}

\smallskip

{\it Proof of Claim 1.} Blowing up $\PP^{n+1}$ with center at $P$
yields a fiber bundle $\xi:\widehat{\PP}^{n+1}\to\PP^n$ with fiber
$\PP^1$. We let $\Sy_{e}(\xi)$ denote the $e$th symmetric
power\footnote{That is the $e$th Cartesian power factorized by the
natural action of the symmetric group of degree $e$.} of $\xi$
over $\PP^n$. Its fiber over a point $Q\in\PP^n$ consists of all
effective divisors on $\xi^{-1}(Q)\cong\PP^1$ of degree $e$. Given
a partition
$$e=\sum_{i=1}^k n_i\qquad\mbox{ with}\qquad
1\le n_1\le n_2\le\ldots\le n_s\,$$ we let $\Sigma_{\bar n}$,
where $\bar n=(n_1,\ldots,n_s)$, denote the closed subbundle of
$\Sy_{e}(\xi)$ whose fiber over $Q$ consists of all effective
divisors on $\xi^{-1}(Q)$ of the form
$$\sum_{i=1}^s n_i[p_i],\qquad\mbox{where}\quad
p_i\in\xi^{-1}(Q)\,.$$ We also let
$$\Sigma_k=\bigcup_{\bar n: n_k\ge 2}
\Sigma_{\bar n}\,.$$ The restriction map
$$\rho:f\longmapsto f^*(0)\cdot (PQ),\qquad Q\in Y\,,$$
associates to $f$ a section $\rho(f)$ of $\Sy_{e}(\xi)$ over $Y$.
It is easily seen that $f\in\F(n+1,e)$ belongs to $\F(Y,e,k)$ if
and only if $\rho (f)$ meets $\Sigma_k$.

We claim that the set, say, $\Gamma_{e,k}$ of all sections of
$\Sy_{e}(\xi)|_Y$ meeting $\Sigma_k$ is a Zariski closed subset of
$\Gamma (Y, {\mathcal O}(\Sy_{e}(\xi)|_Y))$. More generally, given
projective varieties $X$ and $Y$ and a subvariety $S\subset Y$,
the set ${\mathcal M}_S$ of all morphisms $f:X\to Y$ such that the
image $f(X)$ meets $S$ is a Zariski closed subset of $\Mor (X,Y)$.
Indeed, let us consider the incidence relation
$$I=\{(f,x,y)\in \Mor (X,Y)\times X\times Y\,|\,f(x)=y\}\,.$$
Then ${\mathcal M}_S=\pi_1(\pi_3^{-1}(S)\cap I)$ is Zariski
closed, as claimed.

Consequently, $\PP(Y,e,k)$ is Zariski closed in $\PP(n+1,e)$, as
stated. \qed

\smallskip

{\it Claim 2. $\PP(n+1,e)\setminus \PP(Y,e,k)\neq \emptyset$ if
$n=3$.}

\smallskip

Indeed, it is easy to see that the union $X'$ of $e$ planes in
$\PP^3$ in general position belongs to this complement. \qed

\smallskip

Presumably the same holds  in higher dimensions for unions of $e$
hyperplanes in general position. However the latter is much less
evident, so we've chosen above a different approach.\erem

\section{Examples}
\bthm\label{MT} Let $Y_0$ be a Kobayashi hyperbolic hypersurface
in $\PP^n$ ($n\ge 2$), where $\PP^n$ is realized as the hyperplane
$H=\{z_{n+1}=0\}$ in $\PP^{n+1}$. Then a general small deformation
$X_\varepsilon\subseteq\PP^{n+1}$ of the double cone
$X_0=2{\widehat Y}_0$ over $Y_0$ is Kobayashi hyperbolic. \ethm

\bproof Suppose the contrary. Then letting $X_\infty$ be a general
hypersurface of degree $2d=2\deg Y_0$ and $(X_t)_{t\in\PP^1}$ the
pencil generated by $X_0$ and $X_\infty$, we can find a sequence
$\varepsilon_n\longrightarrow 0$ and a sequence of Brody curves
$\varphi_n:\C\to X_{\varepsilon_n}$ such that
$\varphi_n\longrightarrow \varphi$, where $\varphi:\C\to {\widehat
Y}_0$ is non-constant. We let $\pi:{\widehat Y}_0\dashrightarrow
Y_0$ be the cone projection. Since $Y_0$ is assumed to be
hyperbolic we have $\pi\circ\varphi=\,\,$cst. In other words
$\varphi (\C)\subseteq l$, where $l\cong\PP^1$ is a generator of
the cone ${\widehat Y}_0$.

Letting $Y_0=f_0^*(0)$, where $f$ is a homogeneous form of degree
$d$ in $z_0,\ldots,z_n$, we note that $\bigtriangledown f_0^2
|_{{\widehat Y}_0}=0$. If $l$ and $X_\infty$ meet transversally in
a point $\varphi (t)\in l\cap X_\infty$ then $d\varphi(t)=0$ by
virtue of Lemma \ref{prop}.

Since $Y_0\subseteq \PP^n$ is hyperbolic and $n\ge 2$ we have
$d\ge n+2$. In particular $$\deg X_\infty=2d\ge 2n+4\ge 2\dim Y_0
+5\,.$$ By Proposition \ref{trans}, $l$ and $X_\infty$ meet
transversally in at least 5 points. Hence the nonconstant
meromorphic function $\varphi:\C\to l\cong\PP^1$ possesses at
least 5 multiple values. Since the defect of a multiple value is
$\ge 1/2$, this contradicts the Defect Relation.
\end{proof}

\brem\label{rem} Given a hyperbolic hypersurface $Y\subseteq
\PP^n$ of degree $d$, Theorem \ref{MT}  provides a hyperbolic
hypersurface $X\subseteq \PP^{n+1}$ of degree $2d$. Iterating the
construction yields hyperbolic hypersurfaces in $\PP^n$ $\forall
n\ge 3$. However, their degrees $d(n)$ grow exponentially with
$n$, whereas the best asymptotic achieved so far is
$d(n)=4(n-1)^2$ (see e.g., \cite{SZ2}). \erem

\bexa\label{exa1} Let $C\subseteq \PP^2$ be a hyperbolic curve of
degree $d\ge 4$, and let ${\widehat C}\subseteq \PP^3$ be a cone
over $C$. Then a general small deformation of the double cone
$X_0=2{\widehat C}$ is a Kobayashi hyperbolic surface in $\PP^3$
of even degree $2d\ge 8$. \eexa

 The following degeneration principle can be proved along
the same lines as Theorem \ref{MT}.

\bprop\label{prop2} Let $(X_t)_{t\in\PP^1}$ be a pencil of
hypersurfaces in $\PP^{n+1}$ generated by two hypersurfaces $X_0$
and $X_\infty$ of the same degree $d$, where $X_0=kQ$ with $k\ge
2$ for some hypersurface $Q\subseteq \PP^{n+1}$, and
$X_\infty=\bigcup_{i=1}^{d} H_{a_i}$ ($a_1,\ldots,a_d\in\PP^1$) is
the union of $d\ge 5$ distinct hyperplanes from a pencil of
hyperplanes $(H_a)_{a\in\PP^1}$. If a sequence of entire curves
$\varphi_n: \C\to X_{\varepsilon_n}$, where $\varepsilon_n\to 0$,
converges to an entire curve $\varphi: \C\to X_0$, then
$\varphi(\C)\subseteq X_0\cap H_a$ for some $a\in\PP^1\,.$ \eprop

\bexas\label{exs} Given a pencil of planes $(H_a)_{a\in\PP^1}$ in
$\PP^3$, using Proposition \ref{prop2} one can deform
\begin{enumerate}\item[$\bullet$] $X_0=5Q$, where $Q\subseteq\PP^3$ is a
plane, \item[$\bullet$]  a triple quadric $X_0=3Q\subseteq\PP^3$,
or \item[$\bullet$]  a double cubic, quartic, etc.
$X_0=2Q\subseteq\PP^3$ \end{enumerate} to an irreducible surface
$X_\varepsilon\in\langle X_0,X_\infty\rangle$ of the same degree
$d$, where as before $X_\infty=\bigcup_{i=1}^{d} H_{a_i}$, so that
every limiting entire curve $\varphi:\C\to X_0$ is contained in a
section $X_0\cap H_a$ for some $a\in\PP^1$. \eexas

The famous Bogomolov-Green-Griffiths-Lang Conjecture on strong
algebraic degeneracy (see e.g., \cite{BO,GG}) suggests that every
surface $S$ of general type possesses only finite number of
rational and elliptic curves and, moreover, the image of any
nonconstant entire curve $\varphi:\C\to S$ is contained in one of
them. In particular, this should hold for any smooth surface
$S\subseteq\PP^3$ of degree $\ge 5$, which fits the Kobayashi
Conjecture. Indeed, by Clemens-Xu-Voisin's Theorem, a general
smooth surface $S\subseteq\PP^3$ of degree $\ge 5$ does not
contain rational or elliptic curves, hence it should be hyperbolic
provided that the above conjecture holds indeed.

Anyhow, the  deformation method leads to the following result,
which is an immediate consequence of Proposition \ref{prop2}.

\bcor\label{coro} Let $S\subseteq\PP^3$ be a surface and $Z\subset
S$ be a curve such that the image of any nonconstant entire curve
$\varphi:\C\to S$ is contained in $Z$ \footnote{The latter holds,
for instance, if $S$ is hyperbolic modulo $Z$.}. Let $X_\infty$ be
the union of $d=2\deg S$ planes from a general pencil of planes in
$\PP^3$. Then any small enough linear deformation $X_\varepsilon$
of $X_0=2S$ in direction of $X_\infty$ is hyperbolic. \ecor

Along the same lines, Proposition \ref{prop2} can be applied in
the following setting.

\bexa\label{ex3} Let us take for $X_0$ a double cone in $\PP^3$
over a plane hyperbolic curve of degree $d\ge 4$, and  for
$X_\infty$ the union of $2d$ distinct planes from a general pencil
of planes $(H_a)_{a\in\PP^1}$. Then small deformations
$X_\varepsilon$ of $X_0$ in direction of $X_\infty$ provide
examples of hyperbolic surfaces of any even degree $2d\ge 8$. In
suitable coordinates in $\PP^3$ such a surface can be given by
equation
\begin{equation}\label{1}
Q(X_0,X_1,X_2)^2-P(X_2,X_3)=0\,,\end{equation} where $P,\,Q$ are
generic homogeneous formes of degree $d=2k$ and $k$, respectively.
The latter are actually the Duval-Fujimoto examples \cite{Du2,
Fu}.\eexa

Let us finally turn to the Kobayashi problem on hyperbolicity of
complements of general hypersurfaces. By virtue of
Kiernan-Kobayashi-M.\ Green's version of Borel's Lemma, the
complement $\PP^n\setminus L$ of the union $L=\bigcup_{i=1}^{2n+1}
L_i$ of ${2n+1}$ hyperplanes in $\PP^n$ in general position is
Kobayashi hyperbolic. In particular, this applies to the union $l$
of 5 lines in $\PP^2$ in general position. Moreover \cite{Za1} $l$
can be deformed to a smooth quintic curve with hyperbolic
complement via a small deformation. This deformation proceeds in 5
steps and neither is linear nor very generic. So the following
question arises.

\bsit\label{q3} {\bf Question.}   Let $L$ ($M$) stands for the
union of $2n+1$ ($2n-1$, respectively) hyperplanes in $\PP^n$ in
general position. Is the complement of a general small linear
deformation of $L$ Kobayashi hyperbolic? Is a general small linear
deformation of $M$ Kobayashi hyperbolic? In particular, does the
union of 5 lines in $\PP^2$ (of 5 planes in  $\PP^3$) in general
position admit a general small linear deformation to an
irreducible quintic curve with hyperbolic complement (to a
hyperbolic quintic surface, respectively)? \esit

\end{document}